\documentclass[11pt,a4paper,reqno]{amsart} 

\usepackage{amssymb,graphicx}

\usepackage{amssymb}

\newtheorem{theorem}{Theorem}[section]

\newtheorem{proposition}[theorem]{Proposition}

\newcommand{\R}{\mathbb{R}}
\newcommand{\RP}{\mathbb{RP}}
\newcommand{\PP}{\mathbb{P}}
\newcommand{\om}{\omega}

\newcommand{\ad}{\mathrm{ad}}
\newcommand{\rk}{\mathrm{rk}}

\newcommand{\const}{\mathrm{const}}
\newcommand{\spn}{\mathrm{span}}

\newcommand{\XX}{\mathcal{X}}

\newcommand{\VV}{\mathcal{V}}
\newcommand{\DD}{\mathcal{D}}
\newcommand{\II}{\mathcal{I}}
\newcommand{\JJ}{\mathcal{J}}
\newcommand{\tit}{\tilde t}
\newcommand{\p}{\partial}
\newcommand{\be}{\begin{equation}}
\newcommand{\ee}{\end{equation}}

\begin{document}\date{October 15, 2013}
\title{Point invariants of third-order ODEs and hyper-CR Einstein-Weyl structures}

\author{Maciej Dunajski}
\address{Department of Applied Mathematics and Theoretical Physics\\ 
University of Cambridge\\ Wilberforce Road, Cambridge CB3 0WA, UK.
}\email{m.dunajski@damtp.cam.ac.uk}
\author{Wojciech Kry\'nski}
\address{Institute of Mathematics of the  Polish Academy of Sciences\\
\'Sniadeckich 8, 00-956 Warsaw, Poland.}
\email{krynski@impan.pl}

\begin{abstract}
We characterize Lorentzian three--dimensional hyper-CR Einstein-Weyl structures in terms of invariants of the associated third order ordinary differential equations.
\end{abstract}
\maketitle

\section{Introduction}
Given a third-order ODE for $x=x(t)$
$$
x'''=F(t,x,x',x'')\eqno{(F)}
$$
the W\"unschmann invariant $W$ (introduced in \cite{Wun}, and given by
formula (\ref{wun}) below) is defined on the space of 2-jets $J^2(\R,\R)$.
The vanishing of $W$  is the necessary and sufficient condition for the three--dimensional solution space $B$ of $(F)$ to be equipped with a conformal structure $[h]$ of Lorentzian signature
and such that the hyper--surfaces $\xi\subset  B$ corresponding to points in a  two--dimensional manifold $Z=J^0({\R, \R})$
are isotropic.
The W\"unschmann  condition $W=0$ is invariant under contact transformations
of $\PP(TZ)$.

The Cartan invariant $C$ (introduced in \cite{C0} and
given by formula (\ref{cartan_inv}) below)
 is only well defined if $W=0$.
Its vanishing guarantees that the hyper--surfaces $\xi$ defined above are also totally geodesic with respect
to some torsion--free connection $D$ on $B$. It then follows automatically that the connection $D$ and the conformal  structure $[h]$ form a Weyl structure: the null geodesics of $[h]$ are also geodesics of $D$.
It also follows from the Frobenius theorem that the pair $(D, [h])$ satisfies the Einstein--Weyl equations:
the trace--free part of the symmetrised Ricci tensor of $D$ vanishes \cite{C}.
The Cartan condition $C=0$ is only invariant under the point transformations of $Z$.

In \cite{DMT} it was argued that three--dimensional 
Einstein--Weyl equations provide the universal
setup for dispersionless integrable systems. They are naturally equipped with dispersionless Lax pair, and, in the real analytic context, there exists a twistor construction for the solutions. A separate
evidence for this has recently been provided using a completely different approach \cite{FK}.

Many dispersionless integrable systems, like the dKP equation \cite{DMT}, 
 the $SU(\infty)$ Toda equation
\cite{Ward}, the hyper--CR equations \cite{D}, and others \cite{Bog, FK} arise from special Einstein--Weyl structures. It is therefore interesting to seek the characterisation of these integrable systems
in terms of additional invariants of the ODE $(F)$. Such characterisation should reduce the
number of arbitrary functions in the general solution: it is known that in the real analytic category a general
Einstein--Weyl structure depends on four arbitrary functions of two variables. On the other hand solutions
to the integrable systems listed above depend only of two such functions.

In this paper we shall characterise the hyper--CR Lorentzian Einstein--Weyl structures in terms of point
invariants of the ODE $(F)$. The hyper--CR structures form a subclass of all Einstein-Weyl structures
characterised by the existence of the Lax representations without derivatives w.r.t. the spectral parameter. Equivalently, at least in the real analytic category,  they are characterised by the existence of a holomorphic fibration
of the associated twistor space $Z$ over a projective line \cite{D}.  In the ODE approach this projective line
is coordinatised by the independent variable $t$ in $(F)$. Thus the characterisation in terms of a restricted
transformations preserving $t$ is relatively easy to come by (it is given by vanishing
of $K_0$ and $K_1$ defined by (\ref{eqK})). We shall go a step further, and provide a full
characterisation of the hyper--CR EW condition under point transformations 
(i. e. diffeomorphisms of $Z$).
\begin{theorem}\label{thm2}
A third-order ODE $(F)$  defines  a Lorentzian  hyper-CR Einstein-Weyl structure $(D, [h])$ on its solutions space $B$ such that surfaces $(x=const, t=const)$ in $B$ are isotropic with respect to $[h]$ and totally geodesics with respect to $D$ 
if and only if $(F)$ is either point equivalent to the trivial equation, or
$$
W=0, \quad \II=0, \quad \JJ=0,
$$
where the $\II$ and $\JJ$ are respectively a 3-form and a 2-form
on $J^2(\R, \R)$ which are relative point invariants of $(F)$.
\end{theorem}
The explicit expressions for $\II$ and $\JJ$ will be provided in 
Section \ref{sectionew} (formulae (\ref{three_form}) and (\ref{2_form_j}))-
they involve derivatives of $F$ up to 8th order. The 3-form $\II$ is 
invariant and well defined as long as $(F)$ is not point equivalent to a 
trivial ODE. The 2-form $\JJ$ is invariant and well defined if $\II=0$
and $W=0$. The vanishing of $\JJ$ implies the vanishing of the Cartan 
invariant.

\section{ODEs, distributions and vector fields}
\subsection{Preliminaries}
A distribution $\DD$ on a manifold $M$ is a sub-bundle of the tangent bundle $TM$. Let $\Gamma(\DD)$ denote the set of smooth sections of $\DD$. If $\DD$ and $\tilde \DD$ are two distributions on $M$ then their Lie bracket at $x\in M$ is the subspace of $T_xM$ defined by
$$
[\DD,\tilde \DD](x)=\spn\{[X,\tilde X](x)\ |\ X\in\Gamma(\DD), \tilde X\in\Gamma(\tilde \DD)\},
$$
where $[X,\tilde X]$ is the usual Lie bracket of vector fields. We will work locally in a neighbourhood of $x\in M$ where the function $x\mapsto\rk\,[\DD,\tilde \DD](x)$ is constant so that $[\DD,\tilde \DD]$ is  a distribution which we shall call $\ad_\DD\tilde \DD$. We then inductively 
define $\ad^1_\DD \tilde \DD =\ad_\DD\tilde \DD$ and   $\ad^{i+1}_\DD\tilde \DD=[\DD,\ad^i_\DD\tilde \DD]$. Given a vector field $X$ and a distribution $\DD$ we similarly denote
$\ad_X\DD(x)=[X,\DD](x)=\spn\{[X,Y](x)\ |\ Y\in\Gamma(\DD)\}$ and $\ad^{i+1}_X\DD=[X,\ad^i_X\DD]$.

The  \emph{total derivative} vector field
\be
\label{total_der}
X_F=\partial_t+x_1\partial_{x_0}+x_2\partial_{x_1}+F(t,x_0,x_1,x_2)\partial_{x_2}
\ee
is defined on the space of 2-jets $J^2(\R,\R)$ with  coordinates $(t,x_0,x_1,x_2)$. 
The rank-one distribution on $J^2(\R,\R)$ spanned by $X_F$ given by (\ref{total_der})  will be denoted $\XX_F$, i.e.
$$
\XX_F=\spn\{X_F\}.
$$
The space of 2-jets $J^2(\R,\R)$ is also equipped with the following integrable distributions
\be
\label{distributions_on_j2}
\DD^1=\spn\{\partial_{x_2}\},\qquad\DD^2=\spn\{\partial_{x_1},\partial_{x_2}\},\qquad \DD^3=\spn\{\partial_{x_0},\partial_{x_1},\partial_{x_2}\}
\ee
which are tangent to the fibres of the projections $J^2(\R,\R)\to J^1(\R,\R)$, $J^2(\R,\R)\to J^0(\R,\R)$ and $J^2(\R,\R)\to \R$, respectively.
\vskip 5pt
The distribution $\XX_F$ can alternatively be defined by three 1--forms
\be
\label{framing}
\om_0=dx_0-x_1 dt, \quad \om_1=dx_1-x_2 dt, \quad \om_2=dx_2-Fdt
\ee
annihilating $X_F$ given by (\ref{total_der}). An integral curve
$x=x(t)$ of $(F)$ lifts to a curve
\[
t\rightarrow(t, x_0=x(t), x_1=x'(t), x_2=x''(t))
\]
in $J^2(\R, \R)$, and the one-forms $\om_i, i=0, 1, 2$ vanish on this curve.
We shall extend these 1-forms to a basis of $\Lambda^1(J^2(\R, \R))$ by
$\om_3=dt$.

In what follows we shall consider two classes of transformations
of $(F)$:
\be
\label{point_tr}
\mbox{Point transformations:} \qquad
(t, x)\rightarrow
(\tilde{t}=\tilde{t}(t, x),\; \tilde{x}=\tilde{x}(t, x)).
\ee
\be
\label{rfp}
\mbox{Veronese transformations:}\quad
(t, x)\rightarrow
\Big(\tilde{t}=\frac{at+b}{ct+d}, ad-bc\neq0,\;  \tilde{x}=\tilde{x}(t, x)\Big)
\ee
where $a, b, c, d$ are constant. The transformations 
(\ref{rfp}) form a subclass of fibre-preserving transformations.
They  arise naturally in the
theory of Veronese webs \cite{GZ,DK}, and will play a role
in Section \ref{sectionew}.

 The W\"unshmann invariant is preserved by a larger class
of contact transformations, i. e. maps $f:J^1(\R, \R)\rightarrow 
J^1(\R, \R)$ such that $f^*\omega_0=\lambda\omega_0$ for some function 
$\lambda$. The contact transformations do not preserve
the Cartan invariant and the Einstein-Weyl conditions.
\subsection{Vector fields compatible with point structure}
A section $X$ of $\XX_F$ is the total derivative of an equation contact equivalent to $(F)$ if and only if $\ad_X\DD^1$ and $\ad_X^2\DD^1$ are integrable and $\ad_X^3\DD^1=\ad_X^2\DD^1$ (cf. \cite[Theorem 4.5]{JK}). If in addition $\ad_X^2\DD^1=\DD^2$ then $X$ represents an equation point equivalent to $(F)$. We will say that such a vector field $X$ is \emph{compatible with a point structure} defined by $(\XX_F,\DD^2)$. Clearly, the total derivative $X_F$ defined by equation $(F)$ is compatible with the point structure $(\XX_F,\DD^2)$.

In what follows we will need to make use of the overall scaling freedom in the choice of $X_F$. This is easily determined using a chain rule: Under a point transformation (\ref{point_tr})
the total derivative $X_F$ transforms by 
$
X_{\tilde F}=g^{-1}X_F$, 
where 
\be
\label{prop4}
g={A+Bx_1}
\ee
and  $A$ and $B$ are functions of $t$ and $x_0$ only, given by
$
A=\partial_t\tilde t, B=\partial_{x_0}\tilde t.
$

\subsection{Vector fields compatible with projective structure}

There is another set of distinguished sections of $\XX_F$ defined by the pair $(\XX_F,\DD^1)$. These sections correspond to equations in Laguerre-Forsyth normal form in the linear case. In \cite{K1} they were referred to as \emph{projective vector fields}. A section $X\in\Gamma(\XX_F)$ is a projective vector field if
\begin{equation}\label{c1}
\exists V\neq 0\in\Gamma(\DD^1)\qquad \ad_X^3V=-W V\mod\XX_F,
\end{equation}
for some function $W$. In \cite[Proposition 4.1]{K1} we have proved than given an ODE $(F)$ there always exists a projective vector field $X$.
Moreover, if $g^{-1}X$ is a different section of $\XX_F$ satisfying \eqref{c1} then
\begin{equation}\label{eq4}
3X(g)^2-2gX^2(g)=0.
\end{equation}
In general, the projective vector fields are not related to vector fields compatible with the point structure $(\XX_F,\DD^2)$. In the next Section we will see the interplay between the two notions 
in the case of ODEs corresponding to the Einstein-Weyl structures.

The coefficient $W$ in formula \eqref{c1} is the W\"unschmann invariant 
given by
\be
\label{wun}
W=\p_{x_0}F -\frac{1}{2}X_F(\p_{x_1}F) +\frac{1}{3}\p_{x_1}F\p_{x_2}F+\frac{1}{6}X_F^2(\p_{x_2}F)
 -\frac{1}{3}X_F(\p_{x_2}F)\p_{x_2}F+\frac{2}{27}(\p_{x_2}F)^3.
\ee
If $X\mapsto g^{-1}X$ then $W$ transforms as
$
W\mapsto g^{-3}W,
$
where $g$ is a function satisfying \eqref{eq4}. 
Equation \eqref{eq4} can be interpreted as the Schwartz equation. It follows that on each integral curve of $\XX_F$ there is a canonical projective structure. If the W\"unschmann condition holds, the projective structure descents to the solutions space $B=J^2(\R,\R)/\XX_F$ and defines a Lorentzian conformal metric (a $GL(2)$-structure) on $B$.
\section{Hyper--CR Einstein-Weyl structures and ODEs}
\label{sectionew}

As explained in the Introduction, the seminal result of Cartan \cite{C} is that
the Einstein--Weyl condition on a connection $D$ and a Lorentzian  
conformal class $[h]$ on a three--manifold $B$ is equivalent to the existence
of a two parameter family $Z$ of totally geodesic null surfaces in $B$.
This underlies the Lax pair formulation of the Einstein-Weyl condition
\cite{DMT}: Let $V_1, V_2, V_3$ be independent
vector fields on $B$ such that a contravariant metric in $h\in [h]$ is  
\be
\label{three_metric}
h=V_2\otimes V_2-2(V_1\otimes {V_3}+ V_3\otimes V_1).
\ee
Then there exists a connection $D$ such that $(B, [h], D)$ is Einstein--Weyl
if the rank-2 distribution ${\mathcal{D}}_Z$ on $B\times\RP^1$ spanned by the Lax pair
\be
\label{EW_lax}
L_0=V_1-t V_2+f_0\frac{\p}{\p t},\quad
L_1=V_2-t V_3+f_1\frac{\p}{\p t},
\ee
is Frobenius integrable
for some functions $(f_0, f_1)$ which are cubic polynomials in
$t \in \RP^1$. Conversely, every Einstein--Weyl structure arises from some Lax pair (\ref{EW_lax}). The hyper--CR Einstein Weyl spaces are characterised by the existence
of a Lax pair (\ref{EW_lax}) such that $f_0=f_1=0$. 

In the ODE approach the Lax pair distribution is identified with the distribution 
$\DD^2$ defined in (\ref{distributions_on_j2}). This leads to the double fibration picture
\[
B\longleftarrow {\mathcal F}\longrightarrow Z=J^0(\R, \R),
\]
where the correspondence space ${\mathcal F}$ is $B\times \RP^1$ or equivalently
$J^2(\R, \R)$. 
There are two natural coordinate systems on ${\mathcal F}$. Regarding it as the projective spin bundle of $B$, a point $(b, t)\in {\mathcal{F}}$ corresponds to a point in $b\in B$ lying on a totally
geodesic surface given by $t$. The quotient manifold ${\mathcal F}/{\mathcal{D}}_Z$ is then the twistor space
$Z$
of the Einstein--Weyl structure, defined to be the space of isotropic totally geodesic surfaces.
The space $B$ is the quotient of ${\mathcal{F}}$ by a one dimensional distribution spanned by
$\p/\p t$.
Alternatively, regarding ${\mathcal{F}}$ as the space of 2-jets $J^2(\R, \R)$ with coordinates
$(t, x_0, x_1, x_2)$, the space $Z=J^0(\R, \R)$ arises as the quotient ${\mathcal{F}}/\DD^2$. In this approach
$B$ is a quotient of ${\mathcal{F}}$ by the total derivative $X_F$. Note that the meaning of $\p/\p t$ differs between the two pictures as the partial derivative depends on the choice of remaining coordinates. Nevertheless,
the one--dimensional distributions used to obtain the quotient $B$ in both pictures are equivalent. The ODEs arising from Einstein--Weyl structure have $W=0$, and additionally $C=0$, where
\be
\label{cartan_inv}
C={X_F}^2(\p^2_{x_2} F)-X_F(\p_{x_1}\p_{x_2} F)+\p_{x_0}\p_{x_2} F
\ee
is the Cartan invariant.

In \cite{DK} we have pointed out  that there is a one-to-one correspondence between hyper-CR Einstein-Weyl structures and Veronese webs.
A Veronese web \cite{GZ} on a three--dimensional manifold  $B$ is a one-parameter family of foliations of $B$ by surfaces, such that the normal vector fields to these surfaces 
form a Veronese curve in $\PP (T^* B)$.
The curve in $\PP(T^* B)$ defines the Veronese curve $t\rightarrow V(t)$
in $\PP(TB)$ and the vector field $V(t)$ given by  
\[
t \longrightarrow V(t)=V_1-2t V_2+t^2 V_3,\quad\mbox{where}\quad t\in\RP^1
\]
is isotropic w.r.t. the conformal structure (\ref{three_metric}) for any value of $t$.
The vector fields  $L_0$ and $L_1$ given by (\ref{EW_lax}) with $f_0=f_1=0$ form an orthogonal complement of $V(t)$.
Therefore, for each $t\in \RP^1$ they span a null surface in $\zeta\subset B$ which is totally geodesic with respect to $D$. 
This, as explained in the Introduction, is equivalent to the Einstein--Weyl condition, but now with an extra restriction $f_0=f_1=0$.
Thus the  necessary and sufficient condition for an Einstein--Weyl structure to correspond to a Veronese web is that the associated  Lax pair does not contain derivatives w.r.t. $t$.

It is known, \cite{K2}, that  the Veronese webs correspond to equivalence classes of third order ODEs defined modulo
Veronese transformations  (\ref{rfp}). These equivalence classes are characterised  (albeit not invariantly, in a sense described below) 
by two differential conditions $K_0=K_1=0$. Here, for a given ODE $(F)$,  the functions $K_0$ and $K_1$ are implicitly given by 
\cite[Proposition 2.2]{JK}
\begin{equation}\label{c2}
\exists V\neq  0\in\Gamma(\DD^1)\qquad \ad_{X_F}^3V=-K_0V+K_1\ad_{X_F}V\mod\XX_F.
\end{equation}
In \cite[Proposition 2.2]{JK} it has been shown that this definition 
does not depend on the particular choice of $V$ such
that $\ad_{X_F}^3V=0\mod \DD^2$ . In what follows, we will need the explicit form of  $K_0$ and $K_1$
\begin{eqnarray}\label{eqK}
K_0&=&\partial_{x_0}F-X_F(\partial_{x_1}F) +\frac{1}{3}\partial_{x_1}F\partial_{x_2}F+ \frac{2}{3}X_F^2(\partial_{x_2}F)\nonumber\\
 &&-\frac{2}{3}X_F(\partial_{x_2}F)\partial_{x_2}F+ \frac{2}{27}(\partial_{x_2}F)^3,\nonumber\\
K_1&=&\partial_{x_1}F-X_F(\partial_{x_2}F)+\frac{1}{3}(\partial_{x_2}F)^2.
\end{eqnarray}
These functions already appear in the work of Chern \cite{chern}
(where they are called $Q$ on page 266).
They are related to the W\"unschmann and Cartan invariants by
\begin{equation}\label{eqWC}
W=K_0+\frac{1}{2}X_F(K_1),\qquad C=\frac{3}{2}\p_{x_1}K_1+\p_{x_2}F\p_{x_2}K_1+\frac{3}{2}\p_{x_2}K_0,
\end{equation}
but unlike $W$ and $C$, the expressions $K_1$ and $K_2$ are only invariant under a more restricted class of Veronese transformations
(\ref{rfp}).
Therefore the problem considered in the present paper is reduced to characterisation of point equivalent classes of third-order ODEs containing as a representative an equation with vanishing $K_0$ and $K_1$. 

Let us notice that if $K_1=0$, then $K_0$ coincides with the 
W\"unschmann invariant $W$. In this case $X_F$ is a projective vector field for the pair $(\XX_F,\DD^1)$. Thus to solve our main problem, given an equation
$(F)$, we need to  check if there is a point equivalent $(\tilde F)$ such that the invariants $K_0$ and $K_1$ vanish for $(\tilde F)$. In fact, we can assume from the start  that the W\"unschmann condition $W=0$ holds, and 
look for $(\tilde F)$ such that $K_1=0$. The transformation rule for $K_1$ is as follows (see \cite{K1}): if $X_F\mapsto g^{-1}X_F$ then
\begin{equation}\label{eqK1}
K_1\mapsto g^{-2}K_1 +2g^{-3}{X_F}^2(g)-3g^{-4}{X_F}(g)^2.
\end{equation}
On the right hand side one can recognise the Schwartz derivative from \eqref{eq4} thus it follows that the M\"obius transformations of the independent variable $t$ do not affect $K_1$.

To sum up, we start with an equation $(F)$ which satisfies $W=0$. Then we are looking for a function
$g$ of the form (\ref{prop4})
such that $K_1$ vanishes for $(\tilde F)$, where $X_{\tilde F}=g^{-1}X_F$. But then $X_{\tilde F}$ is also a projective vector field. In this way we proved
\begin{proposition}\label{thm3}
If $(F)$ is a third-order ODE with vanishing W\"unschmann invariant $W$ then it defines a 
hyper-CR Einstein-Weyl structure on its solution space if and only if there exists a vector field $X\in\Gamma(\XX_F)$ which is simultaneously compatible 
with the point structure and the projective structure. 
\end{proposition}
In view of Proposition \ref{thm3}, and formula (\ref{eqK1}) the solution space of $(F)$ admits a hyper--CR Einstein-Weyl structure iff
there exists a non-vanishing function $g$ of the form (\ref{prop4}) which satisfies 
\be
\label{eq5}
-2gX_F^2(g)+3X_F(g)^2=g^2K_1.
\ee
\subsection{Proof of Theorem \ref{thm2}}
Our main Theorem \ref{thm2} is a direct application of
 Proposition \ref{thm3} and  Propositions \ref{propI}, and \ref{propJ},
which we shall establish in this Section. 
We shall split our construction into two steps
\subsubsection*{Step 1 -- necessary conditions}
Substituting the coordinate expression (\ref{total_der}) of the total 
derivative  $X_F$ in \eqref{eq5} and differentiating three times w.r.t.
$x_2$ gives 
\begin{equation}\label{eq6'}
\partial_{x_1}g=-\Psi g,
\end{equation}
where
\begin{equation}\label{eqpsi}
\Psi=\frac{\partial_{x_2}^3K_1}{2\partial_{x_2}^3F}.
\end{equation}
The denominator $\partial_{x_2}^3F$ in $\Psi$ is a point invariant of $(F)$. If it vanishes then formula \eqref{eqK} implies that $\partial_{x_2}^3K_1=0$. In this case $(F)$ is point equivalent to the trivial equation \cite{T}. Therefore, we will assume that $\partial_{x_2}^3F\neq 0$.

If we use \eqref{prop4} to rewrite equation \eqref{eq6'} in terms of $\tit$, we find that $\tit$ is constant along
a vector field $\Psi\p_t+(1+x_1\Psi)\p_{x_0}$. Equivalently, $\tit$ is constant along a rank-3 distribution
\be
\label{rank_3_dist}
\VV=\spn\{V_F,\partial_{x_1},\partial_{x_2}\},\quad \mbox{where}\quad
V_F=\partial_{x_0}+\Psi X_F
\ee
on $J^2(\R, \R)$.  This distribution must therefore be integrable, as it is of co-dimension one.
In what follows, we shall use the framing (\ref{framing}) on 
$\Lambda^1(J^2(\R, \R))$.
Let
\be
\label{alpha_form}
\alpha_F=\om_3-\Psi\om_0
\ee
be a one-form which annihilates $\VV$. Then 
\begin{eqnarray}
\label{three_form}
\II&=&d\alpha_F\wedge\alpha_F\\
&=& -(I_1\om_1+I_2\om_2)\wedge \om_0\wedge\alpha_F,\nonumber
\end{eqnarray}
where 
\begin{equation}\label{defI}
I_1=\partial_{x_1}\Psi-\Psi^2,\qquad 
I_2=\partial_{x_2}\Psi,
\end{equation}
where in the derivation we have repeatedly used (\ref{eq6'}).
Now  the Frobenius theorem implies that $[\VV, \VV]\subset \VV$ iff $I_1=I_2=0$.
The two quantities $I_1$ and $I_2$ are not point invariant, but the  3-form $\II$ is.
To see it we explicitly compute the transformation rule 
$
\alpha_F\mapsto g\alpha_F,
$
which implies that $\II\mapsto g^{2}\II$. We have proved
\begin{proposition}\label{propI}
The 3-form $\II$ given by (\ref{three_form}) is a relative invariant of $(F)$. 
Under a  point transformation (\ref{point_tr})  with $g$ given by (\ref{prop4}),
the 3-form $\II$ transforms as
\be
\label{trans_ii}
\II\mapsto g^{2}\II.
\ee
If a third order ODE defines a hyper-CR structure on its solution space then either it is point equivalent to the trivial ODE, or $\II=0$.
\end{proposition}

We have therefore established that
if $(F)$ defines a hyper-CR Einstein-Weyl structure on its solution space then any $(\tilde F)$, which is point equivalent to $(F)$ and such that $K_1$ vanishes for $(\tilde F)$, has an independent variable $\tit$ such that the foliation $\tit=\const$ is tangent to the distribution $\VV$. In particular,  there exist coordinates $(\tit,\tilde x_0,\tilde x_1,\tilde x_2)$ on $J^2(\R,\R)$, adapted to $(\tilde F)$, such that $\VV=\spn\{\partial_{\tilde x_0},\partial_{\tilde x_1},\partial_{\tilde x_2}\}$. 
\subsubsection*{Step 2 -- sufficient conditions}
We now assume that $(F)$ is an equation such that $\II=0$ and $(\tilde F)$ is a point-equivalent equation such that the foliation $\tilde t=\const$ is tangent to the distribution $\VV$. We will denote the invariant $K_1$ of $(F)$ by $K_1(F)$. So, $K_1(\tilde F)$ is the invariant $K_1$ of $(\tilde F)$. We still do not know if $K_1(\tilde F)=0$. However, we can iterate our procedure once more and try to solve \eqref{eq5} again. The current situation is simpler. Indeed, it follows from our discussion in {\bf Step 1} that we can consider only these transformations of the independent variable that preserve $\VV$, i. e. transformations of the form
$
\hat t=\hat t(\tit).
$
Hence, we are looking for the solutions to \eqref{eq5} in the form
$
g(\tit)={A(\tit)},
$
where $A=\partial_{\tit}\hat t$ for some function $\hat t$. The LHS of  \eqref{eq5} now only depends on $\tilde{t}$. 
Thus, the solution ${g}={g}(\tit)$ exists if and only if the right hand side is a function of 
$\tit$. It means that $K_1(\tilde F)$ is constant on leaves of $\VV$, i. e. that\be
\label{intem_J}
dK_1(\tilde F)\wedge\alpha_{\tilde F}=0,
\ee
where $\alpha_{\tilde F}={g}\alpha_F$ is the annihilator of $\VV$ given by a multiple of (\ref{alpha_form}). To rewrite (\ref{intem_J}) as an invariant condition we need to make it independent on ${g}$.
We therefore eliminate from (\ref{intem_J})  the derivatives of ${g}$ along the sections of (\ref{rank_3_dist}) as follows: $\p_{x_2}{g}=0$, as ${g}$ is given by a special case of (\ref{prop4}).
The expression for $\p_{x_1}{g}$ is given by (\ref{eq6'}). Finally, to evaluate $V_F({g})$ we apply $V_F$ to both sides of (\ref{eq6'}), and use the integrability of the distribution $\VV$ in the form
$[V_F, \p_{x_1}]\in\mbox{span}{(\p_{x_1}, \p_{x_2})}$ to simplify the LHS. This removes all derivatives of ${g}$ from (\ref{intem_J}) leaving 
only a  power of ${g}$ as an overall multiple
of the LHS of (\ref{intem_J}). The condition  (\ref{intem_J}) has only been defined up to an overall multiple, so it is legitimate to rescale it and remove the ${g}$--dependence completely.
This procedure leads to an invariant two--form, which is best written down by expanding the coefficients of $dK_1(\tilde F)$ in  (\ref{intem_J}) in terms of the invariant one--forms associated with the ODE $(F)$.
The resulting expression is
\be
\label{2_form_j}
\JJ=(J_0\om_0+J_1\om_1+J_2\om_2)\wedge\alpha_F,
\ee
where
\begin{eqnarray}
J_0&=&\partial_{x_0}K_1+2X_F(\Psi)K_1+\Psi X_F(K_1)-2X_F^3(\Psi)-2\partial_{x_0}F\Psi,\nonumber\\
J_1&=&\partial_{x_1}K_1+2\Psi K_1-6X_F^2(\Psi)-2\partial_{x_1}F\Psi,\label{defJ}\\
J_2&=&\partial_{x_2}K_1-6X_F(\Psi)-2\partial_{x_2}F\Psi\nonumber.
\end{eqnarray}
Moreover we verify explicitly, that a point transformation $X_F\mapsto g^{-1}X_F$ with $g$ given by (\ref{prop4}) yields $\JJ\mapsto g^{-1}\JJ$. We have established the following result

\begin{proposition}\label{propJ}
Let $(F)$ be third-order ODE such that $W=0$ and $\II=0$, where $\II$ is given by (\ref{three_form}). Then the 2-form $\JJ$ defined by (\ref{2_form_j}) is a relative invariant of $(F)$. Under a point 
transformation  (\ref{point_tr}) with $g$ given by $(\ref{prop4})$ it transforms as
\[
\JJ\mapsto g^{-1}\JJ.
\]
Vanishing of $\JJ$ is a sufficient condition for a solution space
to admit a hyper-CR Einstein--Weyl structure.
\end{proposition}
We have reached the end of the road. Given that $\II$ defined by 
(\ref{three_form}) vanishes, there exists a point transformation such that 
the level sets of the  independent variable in ODE $(F)$ are tangent to the distribution $\VV$. If $\JJ$ also vanishes, then $(F)$ is point equivalent to an ODE
with $K_1=0$. But this, given that the W\"unshmann condition holds,
implies that the solution space $B$ admits a hyper-CR Einstein--Weyl structure.
Thus, as both  $\II$, and $\JJ$ are relative point invariants of $(F)$,
their vanishing (together with $W$) is sufficient and necessary
for the statement of Theorem \ref{thm2} to hold. 

Finally, we point out that if $W=\II=\JJ=0$, then the Cartan invariant (\ref{cartan_inv}) also vanishes. This is because a direct computation using 
(\ref{defJ}) gives
\[
C=\frac{3}{2}\partial_{x_2}W-\frac{3}{4}X_F(J_2)+\frac{3}{4}J_1+\frac{1}{4}\partial_{x_2}FJ_2.
\]
\section{Examples}

{\bf 1.} The following equation defines Einstein-Weyl geometry which is not of hyper-CR type 
\[
x'''=\displaystyle 24\,{\frac {{x''}^{3}}{ \left( -3+ \sqrt{9-2\,x'x''} \right) ^{3}}}+12\,{\frac {x'{x''}^{4}}{ \left( -3+ \sqrt{9-2\,x'x''} \right) ^{4}}}.
\]
Indeed, one can verify that
$
W=C=0
$
but $\II\neq 0$. This Einstein--Weyl structure
belongs to the  dispersionless KP (dKP)
class \cite{DMT} which is characterised by the existence
of a conformally weighted vector which is parallel w. r. t. the Weyl 
connection. Characterising the dKP class by point invariants of the ODE is an interesting open problem.
\vskip 2ex
\noindent
{\bf 2.} The equation
$$
x'''=\left(x''\right)^{3/2}
$$
defines  Einstein-Weyl geometry of hyper-CR type. We verify both
$K_1$ and $W$ vanish. Therefore $\II=\JJ=0$. The resulting Einsten--Weyl 
structure is called ${\bf Nil}$. Its conformal class contains a left--invariant
metric on the Heisenberg group \cite{T, DT, D}.
\subsubsection*{Acknowledgements} The work of Wojciech Kry\'nski has been partially supported by the Polish National Science Centre grant ST1/03902.

\end{document}